\documentclass[]{smfart}

\usepackage{amsmath} 
\usepackage{amssymb} 
\usepackage{amsthm} 
\usepackage[english]{babel} 
\usepackage[font=small,justification=centering]{caption} 
\usepackage[nodayofweek]{datetime}
\usepackage[shortlabels]{enumitem} 
\usepackage[T1]{fontenc} 
\usepackage[a4paper]{geometry} 
\usepackage[utf8]{inputenc} 
\usepackage{ifthen} 
\usepackage{mathabx} 
\usepackage{mathtools} 
\usepackage[dvipsnames]{xcolor} 
\usepackage{setspace} 
\usepackage{tikz-cd} 
\usepackage{xfrac} 

\usepackage{amscd}
\usepackage{cleveref}
\usepackage[all]{xy}

\makeatletter
\def\@tocline#1#2#3#4#5#6#7{\relax
  \ifnum #1>\c@tocdepth 
  \else
    \par \addpenalty\@secpenalty\addvspace{#2}%
    \begingroup \hyphenpenalty\@M
    \@ifempty{#4}{%
      \@tempdima\csname r@tocindent\number#1\endcsname\relax
    }{%
      \@tempdima#4\relax
    }%
    \parindent\z@ \leftskip#3\relax \advance\leftskip\@tempdima\relax
    \rightskip\@pnumwidth plus4em \parfillskip-\@pnumwidth
    #5\leavevmode\hskip-\@tempdima
      \ifcase #1
       \or\or \hskip 1em \or \hskip 2em \else \hskip 3em \fi%
      #6\nobreak\relax
    \dotfill\hbox to\@pnumwidth{\@tocpagenum{#7}}\par
    \nobreak
    \endgroup
  \fi}
\makeatother

\newtheorem{thm}{Theorem}[section]
\newtheorem{theorem}[thm]{Theorem} \newtheorem{proposition}[thm]{Proposition} 

\newtheorem{corollary}[thm]{Corollary}

\theoremstyle{definition}

\newtheorem{definition}[thm]{Definition}
\newtheorem{remark}[thm]{Remark}
\newtheorem{question}[thm]{Question}
\newtheorem{convention}[thm]{Convention}
\newtheorem{example}[thm]{Example}

\DeclareMathOperator{\Comm}{\mathrm{Comm}}

\DeclareMathOperator{\Ends}{\mathsf{E}}

\DeclareMathOperator{\dist}{\mathsf{dist}}
\DeclareMathOperator{\Hdist}{\mathsf{hdist}}

\newcommand{\calc}{{\mathcal{C}}}

\newcommand{\calf}{{\mathcal{F}}}

\newcommand{\calp}{{\mathcal{P}}}
\newcommand{\calP}{{\mathcal P}}
\newcommand{\calq}{{\mathcal{Q}}}

\newcommand{\calt}{{\mathcal{T}}}

\newcommand{\ZZ}{\mathbb{Z}}

\usepackage{tikz}
\usetikzlibrary{arrows,quotes}
\tikzstyle{blackNode}=[fill=black, draw=black, shape=circle]
\usetikzlibrary{decorations.pathmorphing}
\tikzset{snake it/.style={decorate, decoration=snake}}
\usetikzlibrary{calc}

\usepackage{centernot}

\newcommand{\Z}{\mathbb{Z}}
\renewcommand{\sl}{\mathrm{sl}}

\usepackage[style=alphabetic, citestyle=alphabetic, sorting=nyvt]{biblatex}
\AtEveryBibitem{\ifentrytype{article}{\clearfield{issn}}{}}
\AtEveryBibitem{\ifentrytype{article}{\clearfield{url}}{}}
\AtEveryBibitem{\ifentrytype{incollection}{\clearfield{url}}{}}
\NewBibliographyString{toappear}
\DefineBibliographyStrings{english}{%
  toappear = {to appear},
}

\renewbibmacro*{in:}{%
  \iffieldundef{pubstate}
    {}
    {\printfield{pubstate}%
     \setunit{\addspace}%
     \clearfield{pubstate}}%
  \printtext{%
    \bibstring{in}\intitlepunct}}
\usepackage{csquotes}

\usepackage{tkz-euclide} 
\usetikzlibrary{matrix}

\makeatletter
\let\@wraptoccontribs\wraptoccontribs
\makeatother

\newcounter{commentcounter}

\title[Quasi-isometries of pairs]{A survey on quasi-isometries of pairs: invariants and rigidity}

\author{Sam Hughes}
\address[Sam Hughes]{Mathematical Institute, Andrew Wiles Building, University of Oxford, Oxford, OX2 6GG, UK}
\email{sam.hughes@maths.ox.ac.uk}
\author{Eduardo Mart\'inez-Pedroza}
\address[Eduardo Mart\'inez-Pedroza]{Department of Mathematics and Statistics, Memorial University of Newfoundland, St. John's, NL, Canada}
\email{emartinezped@mun.ca}
\author{Luis Jorge S\'anchez Salda\~na}
\address[Luis Jorge S\'anchez Salda\~na]{Universidad Nacional Aut\'onoma de M\'exico}
\email{luisjorge@ciencias.unam.mx}

\date{30th December 2021}

\begin{document}

\maketitle
\begin{abstract}
This survey studies pairs $(G,\mathcal{P})$ with $G$ a finitely generated group and $\mathcal{P}$ a (finite) collection of subgroups of $G$. We explore the notion of quasi-isometry of such pairs and the notion of a qi-characteristic collection of subgroups. Both notions are abstractions of phenomena that have appeared repeatedly in the work of several people within the quasi-isometric rigidity realm.
\end{abstract}

\section{Introduction}
Recall that a \emph{quasi-isometry} $X\to Y$ between metric spaces is a coarsely Lipschitz map with a coarsely Lipschitz inverse.  Since Gromov introduced the notion of hyperbolic groups and proved that the property is \emph{geometric} \cite{Gromov1987}, that is, a quasi-isometry invariant, one of the central motivating themes within geometric group theory has been that of quasi-isometric rigidity.  This has traditionally come in one of two forms, either showing that some property of finitely generated groups is a geometric property, or by showing that a group $G$ or class $\calc$ of groups is \emph{quasi-isometrically rigid}.  In the later case, this involves showing that every group quasi-isometric to $G$ (resp. a group in $\calc$) is \emph{virtually-isomorphic} to $G$ (resp. a group in $\calc$).  Here two groups $H$ and $K$ are virtually isomorphic if there exists finite normal subgroups $M\trianglelefteq H$ and $M\trianglelefteq K$ and finite index subgroups $H'\leq H/M$ and $K'\leq K/N$ such that $H'\cong K'$.

The first example of quasi-isometric rigidity is essentially due independently to Freudenthal and Hopf where they show that the number of ends of a group $G$ is $2$ if and only if $G$ is virtually cyclic.  Gromov's polynomial growth theorem shows that the class of nilpotent groups is qi-rigid \cite{Gromov1981} and Pansu showed that $\ZZ^n$ is qi-rigid \cite{Pansu1983}. Since then Schwarz \cite{Schwartz1995,Schwartz1996} showed some non-uniform lattices in semi-simple Lie groups are qi-rigid and the general case was completed by Eskin \cite{Eskin1998} (see also \cite{Drutu2000}).  Qi-rigidity of Fuchsian groups follows from work of Tukia \cite{Tukia1988,Tukia1994}, Casson--Jungreis \cite{CJ1994}, and Gabai \cite{Gabai1991}.  That the class of uniform lattices in a given semi-simple Lie group is qi-rigid due to many authors \cite{Tukia1986,Pansu1989,Chow1996,KleinerLeeb1997,EskinFarb1997}.  The area of qi-rigidity is still an area of intense active research; with work on right angled Artin and Coxeter groups \cite{DaniThomas2017,HuangKleiner2018,Huang2018,BoundsXie2020}, mapping class groups \cite{BehrstockKleinerMinksyMosher2012}, and others \cite{FaMo2000,FaMo2002,Taback2000,Xie2006,Wortman2007,FrigerioLafontSisto2015,TaTo2019,ShWo2021}.

More recent research has made considerable efforts to understand the large scale geometry and rigidity of a group $G$ with respect to a subgroup or collection of subgroups $\calp$.  A large part of the theory and motivation to consider such problems has come from the techniques used to study relatively hyperbolic groups.  However, the theory has outgrown its origins and, in the process of doing so, has raised a large number of exciting questions.  This survey aims to both highlight recent notions  that appear in the literature, look at the state of the art whilst bringing it into a common language, and to compile an extensive open problem list.

Quasi-isometric rigidity with respect to a collection of subgroups has essentially taken two forms (relative geometric properties and quasi-isometric rigidity of subgroups) and both can be expressed in the language of \emph{quasi-isometries of pairs}.  We will give the technical definition in \Cref{sec.QIdefn}, but roughly for groups $G$ and $H$ with collections of subgroups $\calp$ and $\calq$ respectively a \emph{quasi-isometry of pairs} $(G,\calp)\to(H,\calq)$ is a quasi-isometry $G\to H$ which shows the large scale geometry of $G/\calp$ is similar to $H/\calq$.  

The two forms of quasi-isometric rigidity with respect to a collection of subgroups can then stated as follows.  Either, showing that certain properties of the pair $(G,\calp)$ are \emph{relatively geometric}, that is, for any quasi-isometry of pairs $(G,\calp)\to (H,\calq)$ the pair $(H,\calq)$ satisfies the same properties as $(G,\calp)$.  Or, that for some fixed collection $\calp$ every quasi-isometry $G\to G$ extends to a quasi-isometry of pairs $(G,\calp)\to(G,\calp)$.

Throughout this survey we will only consider finitely generated groups $G$ with a chosen word metric $\dist_G$, and finite collections of subgroups $\calp$.
In Section~\ref{sec.QIdefn} we will define and contextualise quasi-isometries of pairs, then raise and discuss a number of fundamental questions concerning the notion.  In \Cref{sec.invariants} we will explore relative geometric properties of groups, that is, properties invariant under quasi-isometries of pairs.  We will pay particular attention to geometric properties of elements of the collection, filtered ends, and filling functions.  In \Cref{section:qi:characteristic} we will explore quasi-isometrically characteristic collections which is a very strong form of quasi-isometric rigidity of a collection of subgroups.

\subsection*{Acknowledgements}
This work has received funding from the European Research Council (ERC) under the European Union’s Horizon 2020 research and innovation programme (Grant agreement No. 850930).  The second author acknowledges funding by the Natural Sciences and Engineering Research Council of Canada, NSERC. The third author was supported by DGAPA-UNAM-PAPIIT-IA101221.

\section{Quasi-isometries of pairs: definition, examples and non-examples}\label{sec.QIdefn}
In this section we will give the definition of a quasi-isometry of pairs and highlight various places the notion has appeared in the literature. We will then raise and discuss questions regarding invariants and rigidity of quasi-isometry of pairs.

\begin{definition}[Quasi-isometry of metric pairs]\label{defn:quasi-isometry-pairs}
Let $X$ and $Y$ be metric spaces, let $\mathcal{A}$ and $\mathcal{B}$ be   collections of subspaces of $X$ and $Y$ respectively.  A  quasi-isometry $q\colon X\to Y$ is a \emph{quasi-isometry  of pairs} $q\colon (X,\mathcal{A}) \to (Y,\mathcal{B})$ if there is $M>0$:
\begin{enumerate}
    \item For any $A\in \mathcal{A}$, 
    the set $\{ B\in \mathcal{B} \colon \Hdist_Y(q(A), B) <M \}$ is non-empty.      
    \item For any $B\in \mathcal{B}$, 
    the set $\{ A\in \mathcal{A} \colon \Hdist_Y(q(A), B) <M \}$ is non-empty. 
\end{enumerate}
In this case, if $q\colon X \to Y$ is a $(L,C)$-quasi-isometry, then $q\colon (X, \mathcal{A})\to (Y,\mathcal{B})$ is called a \emph{$(L,C,M)$-quasi-isometry}. If there is a quasi-isometry of pairs  $(X,\mathcal{A}) \to (Y,\mathcal{B})$ we say  that $(X,\mathcal{A})$ and  $(Y,\mathcal{B})$ are \emph{quasi-isometric pairs}.
\end{definition}

The above definition is implicit in the work of Kapovich and Leeb on the quasi-isometry invariance of the geometric decomposition of Haken manifolds~\cite[\S 5.1]{KaLe97}. These ideas have been used in a similar fashion in other works, for example~\cite{DS05, BDM09, MSW11, LFS15}. Quasi-isometries of pairs have recently attracted the attention of other researchers in group theory,  see for example~\cite{BuHr21, HaHr19, GT21,  MaSa21, HuMaSa21}, in \cite{Ge19}  the notion appears implicitly. 

Let $G$ and $H$ be finitely generated groups and let $\calp$ and $\calq$ be collections of subgroups of $G$ and $H$ respectively.  Let $\Hdist_G$ denote the Hausdorff distance between subsets of $G$, and let $G/\calp$ denote the collection of left cosets $gP$ for $g\in G$ and $P\in\calp$.  Viewing $G$ and $H$ as metric spaces with respect to their finite generating sets, the definition of a quasi-isometry of metric pairs takes the following form.

\begin{definition}[Quasi-isometry of pairs]
For constants $L\geq 1$, $C\geq 0$ and $M\geq 0$, an \emph{$(L,C,M)$-quasi-isometry of  pairs}  $q\colon (G, \mathcal{P})\to (H, \mathcal{Q})$ is an $(L,C)$-quasi-isometry $q\colon G \to H$ such that  the relation
\begin{equation}\label{eq:def-qi-relation}
\{ (A, B)\in G/\mathcal{P} \times H/\mathcal{Q} \colon \Hdist_H(q(A), B) <M  \} \end{equation}
satisfies that the projections to $G/\mathcal{P}$ and $H/\mathcal{Q}$ are surjective. 
\end{definition}

\begin{question}\label{question:qi:pairs}
Given a group $G$ and a collection of subgroups $\calp$, and a quasi-isometry $q\colon G\to H$, under what conditions does there exist a collection $\calq$ of subgroups of $H$ such that $q$ extends to a quasi-isometry of pairs $q\colon (G,\calp)\to (H,\calq)$? 
\end{question}

Below we illustrate how this question has been addressed in particular classes of groups,  see \Cref{thm:qi:characteristic}.

There are many quasi-isometry invariants for discrete groups that have  generalisations for pairs, for example relative hyperbolicity, relative Dehn functions, and a number of cohomological properties such as relative (filtered) ends, relative (Bredon) finiteness properties and cohomological dimension, and relative duality groups.  We will discuss many of these (and others) in this survey.

\begin{question}
What properties of a pair $(G,\calp)$ are invariant under a quasi-isometry of pairs?
\end{question}

For the question above there are already positive and negative results. Since we already mentioned some positive answers, let us mention one relative property that is not a quasi-isometry invariant.

In \cite{Alonso1994} Alonso proved that finiteness properties $\mathsf{F}_n$ and $\mathsf{FP}_n$ are quasi-isometry invariants. A relative version of these finiteness properties are provided by Bredon finiteness properties. Given a group $G$ and a family $\calf$ of subgroups, that is $\calf$ closed under conjugation and under taking subgroups, it is possible to define the relative finiteness properties $\calf$-$\mathsf{F}_n$ and $\calf$-$\mathsf{FP}_n$. In our context, given a collection $\calp$ of subgroups of $G$ we can consider the smallest family of subgroups $\calf_\calp$ of $G$ that contains $\calp$. 
In \cite{LN03} Leary and Nucinkis constructed groups $H\subseteq G$ such that $H$ has finite index in $G$ and such that $H$ is of type $\mathsf{F}$ (in particular torsion-free), $G$ has finitely many conjugacy classes of finite subgroups, $G$ is not of type $\calf$-$\mathsf{F}_n$ for $n\geq 1$ and $\calf$ is the family of finite subgroups. In these examples if we consider $\calp$ as a finite set of representatives of finite subgroups of $G$, then $\calf$ is the smallest family that contains $\calf$. It is not difficult to verify that the inclusion map $H\hookrightarrow G$ is a quasi-isometry of pairs $(H,\{1\})\to (G,\calp)$  that does not preserve relative finiteness properties. 

Note however, for certain quasi-isometrically rigid classes of pairs we might expect to have homogeneous Bredon finiteness properties.  For example, relatively hyperbolic groups are of type $\calf$-$\mathsf{F}_\infty$, where $\calf$ is the family generated by the peripheral and finite subgroups of $G$ \cite{MPP2019}.

\subsection{Quasi-isometry of pairs and relatively hyperbolic groups}
Although much of the theory regarding quasi-isometries of pairs is applicable in much more general situations than relatively hyperbolic groups, much of the groundwork is based on insights garnered from the relatively hyperbolicity.  In this section we will examine these connections and suggest possible avenues for extension.

\begin{theorem}\emph{\cite[Theorem~4.1]{BDM09}} Let $G$ and $G'$ be a relatively hyperbolic groups with respect to the finite collections $\calp$ and $\calp'$ respectively. Assume that every $P\in \calp$ and every $P'\in \calp'$ is non-relatively hyperbolic. Then every quasi-isometry $G\to G'$ is a quasi-isometry of pairs $(G,\calp)\to (G',\calp')$.
\end{theorem}


Roughly, JSJ decompositions give a way to understand graph of group splittings of a group $G$ over a class of subgroups $\calp$.  However, there is not enough space in this survey to discuss JSJ decompositions and the beautiful maths they have generated so we defer the interested reader to \cite{GuLe2017} for background.  We do however take the liberty to mention one particularly interesting rigidity result due to Haulmark--Hruska in the context of relatively hyperbolic groups.

\begin{theorem}\emph{\cite[Corollary~1.3]{HaHr19}}
Let $(G,\calp)$ and $(G',\calp')$ be relatively hyperbolic groups with connected Bowditch boundaries. Every quasi-isometry of pairs $(G,\calp)\to (G',\calp')$ induces a vertex-type preserving isomorphism of JSJ-trees $\calt_G\to \calt_{G'}$.
\end{theorem}

In light of this we raise the following very general question.

\begin{question}
When is a given JSJ-splitting of a group $G$ qi-characteristic?
\end{question}

Some work towards this question can be found in the literature.  Most notably is
Panos Papsoglu's work \cite{Pa05}.  There are also results due to Alexander Margolis in the case of right-angled Artin groups \cite{Margolis2019}, and Alexander Taam and Nicholas W.M. Touikan in the case of word hyperbolic groups with a cyclic JSJ decomposition that has only rigid vertex groups \cite{TaTo2019}.

\section{Some results on quasi-isometry invariants for pairs}\label{sec.invariants}
In this section we survey when properties of a pair $(G,\calp)$ are invariant under quasi-isometry.  In \Cref{sec.invariants.subs} we will explain how geometric properties of subgroups in $\calp$ are qi of pairs invariants.  In \Cref{sec.invariants.ends} we discuss work of the second and third author showing that the set of filtered ends is a qi of pairs invariant in analogy with ends being a geometric property.  In \Cref{sec.invariants.filling} we will discuss work of all three authors showing that relative Dehn functions are a qi of pairs invariant in analogy with the Dehn function being a geometric property.  Throughout we will raise a number of questions

\subsection{Geometric properties of subgroups}\label{sec.invariants.subs}
A property of groups which is invariant under quasi-isometry is called \emph{geometric}.  Examples of these include word hyperbolicity, relative hyperbolicity, the growth type of the Dehn function, the finiteness properties $\mathsf{F}_n$ and $\mathsf{FP}_n$, amenability, and the number of ends of the group.  The following proposition demonstrates that geometric properties of subgroups are rigid under a quasi-isometry of pairs.  That is, for any finitely generated $P\in\calp$, there is some subgroup $Q$ of $H$ which satisfies the same geometric properties as $P$.

\begin{proposition}\emph{\cite[Proposition~2.16]{MaSa21}}
Let $(G,\calp)\to(H,\calq)$ be a qi of pairs, $P\in\calp$, $Q\in\calq$ and $h\in H$.  If $P$ is finitely generated and $\Hdist(q(P),hQ)<\infty$, then $Q$ is finitely generated and $P$ and $Q$ are quasi-isometric.
\end{proposition}

One possible strategy to attack the following question may be via Kevin Li's work on bounded Bredon cohomology \cite{Li2021}.

\begin{question}
Let $(G,\calp)\to(H,\calq)$ be a qi of pairs. If $H$ is amenable relative to $\calq$, then when is $G$ amenable relative to $\calp$?
\end{question}

A finitely presented group $G$ is \emph{coherent} if every finitely generated subgroup is finitely presented, otherwise $G$ is \emph{incoherent}.  A long standing folklore problem asks if every finitely presented infinite group with positive Euler characteristic is incoherent (see \cite{Wise2020} for an extensive problem list and \cite[Page 734]{Newman1973}).  The motivation for this is that standard example of an incoherent group is $F_2\times F_2$ where the incoherence is \emph{witnessed} by the kernel $K$ of the map to $\ZZ$ sending every generator to $1$.  The property of being (in)coherent is easily seen to be a commensurability invariant.  However, a closely related open question is whether a group $G$ being (in)coherent is a geometric property.  Since this property concerns subgroups of $G$ a natural avenue to attack this question is via quasi-isometries of pairs:

\begin{question}
Let $G$ be a finitely presented incoherent group with witness $K$.  Does every quasi-isometry $q\colon G\to G$ extend to a quasi-isometry of pairs $q\colon (G,\calp)\to(G,\calp)$ such that $K\in\calp$? Notice that $\calp$ depends on $q$.
\end{question}

We highlight the special case of $F_2\times F_2$ and $K$ defined above.

\begin{question}
Does every quasi-isometry $q\colon F_2\times F_2\to F_2\times F_2$ extend to a quasi-isometry of pairs $q\colon (F_2\times F_2,\calp)\to(F_2\times F_2,\calp)$ such that $K\in\calp$? Notice that $\calp$ depends on $q$.
\end{question}

Here one possible strategy could be to relate the problem to cocompact lattices in a product of two locally finite trees.  Such lattices can be virtually a direct product of free groups, arithmetic lattices in $p$-adic Lie groups, groups with no finite quotients \cite{Wise2007}, or even simple groups \'a la Burger and Mozes \cite{BM97,BM00a,BM00b}.  An example of an incoherent simple group is provided by Rattaggi \cite{Rattaggi2005}.

\subsection{Filtered ends}\label{sec.invariants.ends}



Let $G$ be a finitely generated group and let $P$ be a subgroup.  The  \emph{number of filtered ends $\tilde e(G,P)$} of the pair $(G, P)$ was introduced by Bowditch~\cite{Bow02}, under the name of \emph{coends}, in his study of JSJ splittings of one-ended groups. The number of filtered ends coincides with \emph{the algebraic  number of  ends of the pair $(G,P)$} introduced by  Kropholler and Roller~\cite{KrRo89}, see~\cite{Bow02} for the equivalence. The number of filtered ends does not coincide with the  \emph{number of relative ends $e(G,P)$} introduced by Houghton~\cite{Hou74}, but there are several relations including the inequality $e(G,P) \leq \tilde e(G,P)$ and equality in the case that $P$ is  normal and finitely generated, for an account see Geoghegan's book~\cite[Chapter~14]{Ge08}. Roughly speaking the set of filtered ends $\Ends(X,C)$ for a path-connected metric pair $(X,C)$ is the inverse limit
 $\varprojlim C(X-C^\mu)$, where $Y^\mu$ is the $\mu$-neighborhood of $Y$ and $\mu>0$, and $C(X-C^\mu)$ is the set of connected components of $X-C^\mu$. The number of filtered ends of $(X,Y)$ is the cardinality of $\Ends(X,C)$. It follows that whenever $X$ is a proper metric space and $C$ is compact, then $\Ends(X,C)$ is the usual set of ends of $X$. In \cite{MaSa21} the authors give a more technical definition for filtered ends that deals with general metric pairs.

\begin{proposition}\emph{\cite[Proposition~5.30]{MaSa21}}\label{cor:qiinvariance:relativeends}
Let $X$ and $Y$ be  metric spaces, and $C\subseteq X$ and $D\subseteq Y$. If $f\colon X\to Y$ is a quasi-isometry such that $\Hdist(f(C),D)$ is finite, then $\Ends(f)\colon \Ends(X,C)\to\Ends(Y,D)$ is a bijection.
\end{proposition}

In other words, the above proposition is telling us that the set of filtered ends $\Ends(X,C)\to\Ends(Y,D)$ is quasi-isometry invariant for pairs of the form $(X,\{C\})$.

\subsection{Filling functions}\label{sec.invariants.filling}

For a pair $(G, \calp)$, Osin introduced the notions of \emph{finite relative presentation} and \emph{relative Dehn function} $\Delta_{G,\calp}$ as natural  generalizations of their standard counterparts for finitely generated groups, see~\cite{Osin06}. These notions characterise relatively hyperbolic pairs $(G,\calp)$ as the ones which are relatively finitely presented and have relative Dehn function bounded from above by a linear function. By quasi-isometric rigidity of relative hyperbolicity, among relatively finitely presented pairs,   quasi-isometries of pairs preserve  having linear relative Dehn function.

\begin{convention}[$\Delta_{G,\calp}$  is well-defined]
By \emph{$\Delta_{G,\calp}$  is well-defined} we mean that $G$ is finitely presented relative to $\calp$ and the relative Dehn function $\Delta_{G,\calp}$ takes only finite values with respect to a finite relative presentation of $G$ and $\calp$. From here on, when we refer to a relative Dehn function, we always assume that it has been defined using a finite relative presentation.
\end{convention}

Let $\calp$ be a collection of subgroups of group $G$.  A \emph{refinement} $\calp^\ast$ of $\calp$ is a set of representatives of conjugacy classes of the collection of subgroups $\{\Comm_G(gPg^{-1}) \colon P\in\calp \text{ and } g\in G \}$ where $\Comm_G(P)$ denotes the commensurator of the subgroup $P$ in $G$.

\begin{thm}\label{thmx:DehnQI}
Let $(G,\calp)\to(H, \calq)$ be a quasi-isometry of pairs and let $\calp^*$ be a refinement of $\calp$. If the relative Dehn function $\Delta_{H,\calq}$ is well-defined, then $\Delta_{G,\calp^*}$ is well-defined and $\Delta_{G,\calp^\ast} \asymp \Delta_{H,\calq}$.
\end{thm}

Recall that a simplicial graph $\Gamma$ is fine if, for each natural number $n$, every edge is contained in finitely many circuits of length $n$. In \cite[Theorem~E]{HuMaSa21}, the authors prove that, for a relatively finitely presented pair $(G,\calp)$, the relative Dehn function $\Delta_{G,\calp}$ is well-defined if and only if the coned-off Cayley graph $\hat \Gamma (G,\calp)$ is fine. In a sense this characterisation is given in geometrical terms, and it makes sense to ask the following.

\begin{question}
Is there an algebraic characterization for pairs $(G,\calp)$ for which the relative Dehn function $\Delta_{G,\calp}$ is well-defined? 
\end{question}

The following is an very explicit example of a pair $(G,\calp)$ for which the relative Dehn function is well defined provided by A. Minasyan. Here $BS(k,l)$ stands for the Baumslag-Solitar group given by the presentation $\langle a,t \mid t a^k t^{-1}=a^l\rangle$.

\begin{thm}\emph{\cite[Theorem~ A.1]{HuMaSa21}} \label{thm:BS_rel_DF-character} Let $G=BS(k,l)$, for some non-zero integers $k,l$.
The relative Dehn function $\Delta_{G,\langle t \rangle}$ is well-defined if and only if $k$ does not divide $l$ and $l$ does not divide $k$.
\end{thm}

There is a rich class of  pairs $(G,\calp)$ with well-defined relative Dehn function.
Hyperbolically embedded subgroups were introduced in~\cite{DGOsin2017} by Dahmani, Guirardel and Osin.  Given a group $G$, $X\subset G$ and $H\leq G$, let $H\hookrightarrow_h (G,X)$ denote that $H$ is a hyperbolically embedded subgroup of $G$ with respect to $X$. 


\begin{thm}\label{thmx:last}
Let $G$ be a finitely presented group and $H\leq G$ be a subgroup. If $H\hookrightarrow_h G$ then the relative Dehn function $\Delta_{G,H}$ is well-defined. 
\end{thm}

\section{QI-characteristic collections: definition, examples and non-examples}\label{section:qi:characteristic} 
In this section we will survey quasi-isometrically characteristic collections as introduced by the second two authors in \cite{MaSa21}.  Roughly, a  quasi-isometrically characteristic collection is a collection of subspaces or subgroups satisfying a strong version of quasi-isometric rigidity.  We will now give the technical defintion.

\begin{definition}\label{def:qicharacteristic}
Let $X$ be a metric space with metric $\dist$. A collection of subspaces $\mathcal{A}$ is called \emph{quasi-isometrically characteristic}, or for short \emph{qi-characteristic}, if the following properties hold:
\begin{enumerate}
\item For any $L\geq 1$ and $C\geq 0$ there is $M=M(L,C)>0$ such that any $(L,C)$-quasi-isometry $q\colon X\to X$ is an $(L,C,M)$-quasi-isometry of pairs $q\colon (X, \mathcal{A}) \to (X, \mathcal{A})$.  

\item Every bounded subset $B\subset X$ intersects only finitely many non-coarsely equivalent elements of $\mathcal{A}$; where $A, A' \in \mathcal{A}$ are coarsely equivalent if their Hausdorff distance is finite.

\item  For any $A\in \mathcal{A}$ the set $\{A'\in \mathcal{A} \colon \Hdist(A, A')<\infty\}$ is bounded as a subspace of $(\mathcal{A}, \Hdist)$.
\end{enumerate}
\end{definition} 

Note that the property of being a qi-characteristic collection is invariant under quasi-isometry of pairs, in the sense that given $q\colon (X, \mathcal{A}) \to (Y, \mathcal{B})$ a  quasi-isometry of pairs, 
then $\mathcal{A}$ is qi-characteristic if and only if $\mathcal{B}$ is qi-characteristic.  A number of other elementary  properties were established in \cite{MaSa21}.

\begin{definition}\label{def:group:qicharacteristic} Let $G$ be a finitely generated group, and let $\mathcal{P}$ be a collection of subgroups of $G$. The collection $\mathcal{P}$ is  \emph{qi-characteristic} if $G/\mathcal{P}$ is a qi-characteristic collection of subspaces of $G$. 

A subgroup of $G$ is a \emph{qi-subcharacteristic subgroup} if it belongs to a qi-characteristic collection of subgroups of $G$.
\end{definition}

In this section we list several examples of qi-characteristic collections as well as some recent results. The notion of qi-characteristic subgroup provides a partial positive answer to Question~\ref{question:qi:pairs}.  

\begin{theorem}\emph{\cite[Theorem~1.1]{MaSa21}}\label{thm:qi:characteristic}
Let $G$ be a finitely generated group, let $\mathcal{P}$ be a finite qi-characteristic collection of subgroups of $G$.  If $H$ is a finitely generated group and $q\colon G \to H$ is a quasi-isometry, then  there is a qi-characteristic collection of subgroups $\mathcal{Q}$ of $H$ such that $q\colon (G, \mathcal{P}) \to (H, \mathcal{Q})$ is a quasi-isometry of pairs.
\end{theorem}

The qi-characteristic conditon imposes algebraic constrains on the subgroups of the collection.

\begin{theorem}\emph{\cite[Theorem~2.9]{MaSa21}}\label{thm:qi-characteristic-groups2}
Let $G$ be a finitely generated group. A collection of subgroups $\mathcal{P}$ is qi-characteristic if and only if 
\begin{enumerate}
    \item For any $L\geq 1$ and $C\geq 0$ there is $M=M(L,C)>0$ such that any $(L,C)$-quasi-isometry $q\colon G\to G$ is an $(L,C,M)$-quasi-isometry of pairs $q\colon (G, \mathcal{P}) \to (G, \mathcal{P})$.

    \item $\mathcal{P}$ is finite.
    
    \item Every $P\in \mathcal{P}$ has finite index in its commensurator.
\end{enumerate}
\end{theorem}

An immediate consequence of the previous result is that a finite subgroup of a finitely generated infinite group is not qi-subcharacteristic.

 \subsection{Qi-characteristic collections of fundamental groups of certain manifolds}
A \emph{Haken manifold} is a compact orientable irreducible 3-manifold that contains an orientable incompressible surface.  Such a manifold admits a hierarchy where they can be split up into $3$-balls along incompressible surfaces \cite{Haken1962}.

\begin{theorem}\emph{\cite[Theorem~1.1]{KL97}}
Let $M$ be a Haken manifold with fundamental group $\Gamma$, then the set of fundamental groups of the geometric pieces of $\Gamma$ is a qi-characteristic collection.
\end{theorem}



Following \cite{LFS15}, a \emph{high dimensional graph manifold} is a compact smooth manifold supporting a decomposition into finitely many pieces, heach of which is diffeomorphic to the product of a torus with a finite volume hyperbolic manifold with toric cusps. The pieces of a graph manifolds are glued along the boundary components which happen to be tori that they call \emph{walls}. 

\begin{theorem}\emph{\cite[Proposition~8.35 and Proposition~8.37]{LFS15}}
For a higher graph manifold $M$ with fundamental group $\Gamma$ denote by $\mathcal W$ (resp. $\mathcal P$) the collection of (the embedded) fundamental groups of walls (resp. pieces) of $M$. Then both $\mathcal{W}$ and $\mathcal{P}$ are qi-characteristic collections of $\Gamma$.
\end{theorem}

\subsection{Qi-characteristic collections graph products of finite groups}

Given a simplicial graph $\Gamma$, a subgraph $\Lambda\leq \Gamma$ is \emph{square-complete} if every induced square of $\Gamma$ containing two opposite vertices in $\Lambda$ must be entirely included into $\Lambda$. A \emph{minsquare subgraph} of $\Gamma
$ is a subgraph which is minimal among all the square-complete subgraphs of $\Gamma$ containing at least one induced square.

\begin{theorem}\emph{\cite[Theorem~1.1]{G19}}
Let $\Gamma$ be a finite simplicial graph and $\mathcal{G}$ a collection of finite groups indexed by the vertex set $V(\Gamma)$. Denote by $\Gamma \mathcal{G}$ the corresponding graph product. Then, the collection
\[\{\langle \Lambda  \rangle\leq \Gamma\mathcal{G}\ |\ \Lambda \text{ is a minsquare subgraph of }\Gamma\}\]
where $\langle \Lambda  \rangle$ denotes the subgroup generated by the groups labelling the vertices of $\Lambda$, is a qi-characteristic collection of $\Gamma$ amongst graph products of finite groups.
\end{theorem}

The above theorem applies, in particular, for all right-angled Coxeter groups.

\subsection{Qi-characteristic collections of some wreath products}
Recall that given groups $F$ and $H$. The wreath product $F\wr H$ is by definition the semi-direct product $F^{|H|}\rtimes H$, where $H$ acts on $F^{|H|}$ as follows $h(f_g)_{g\in H}=(f_{hg})_{g\in H}$. These wreath products are also called \emph{lamplighter groups}.

\begin{theorem}\emph{\cite[Theorem 1.18 and Proof of Theorem 7.3]{GT21}}
If $F$ is a finite group and $H$ is a finitely presented one-ended group, then the collection $\{H\}$ is a qi-characteristic collection of $F\wr H$.
\end{theorem}

\subsection{Qi-characteristic collections of relatively hyperbolic groups}
 Following the convention in~\cite{BDM09}, if a group contains no collection of proper subgroups with respect to which is relatively hyperbolic, then we say that the group is \emph{not relatively hyperbolic (NRH)}. The following theorem is a   consequence of a corollary of work by 
Behrstock,   Dru\c{t}u,  and Mosher~\cite[Theorem~4.1]{BDM09} and  \Cref{thm:qi-characteristic-groups2}.  
 
 \begin{theorem}\emph{ \cite[Theorem~3.1]{MaSa21}}\label{cor:BDM-RelHyp}
Let $G$ be a finitely generated group hyperbolic relative to a finite collection $\mathcal{P}$ of NRH finitely generated subgroups. Then $\mathcal{P}$ is a qi-characteristic collection of $G$.
\end{theorem}

\begin{remark}
The NRH hypothesis of \Cref{cor:BDM-RelHyp} is necessary, for instance, if $F$ is a free group of finite rank then a maximal cyclic subgroup is not qi-subcharacteristic. There is a quasi-isometry of $F$ that maps an  infinite geodesic preserved by a non-trivial element of $F$ to a geodesic that is  preserved by no element of $F$.
\end{remark}

\begin{corollary}\label{cor:RemovingQIsubgroup}
Let $G$ and $\mathcal{P}=\{P_0,\ldots, P_n\}$ be as in \Cref{cor:BDM-RelHyp}. Suppose that there is no quasi-isometric embedding $P_i \to P_0$ for $1\leq i\leq n$. Then $\{P_1, \dots , P_n\}$ is a qi-characteristic collection. 
\end{corollary}

\begin{example}
Let $A$ and $B$ be finitely generated NRH groups endowed with word metrics with a common finite subgroup $C$. By \Cref{cor:RemovingQIsubgroup}, if there is  no quasi-isometric embedding $A\to B$, then $\{A\}$ is a qi-characteristic collection of subgroups of $A\ast_C B$.
\end{example}

In contrast to the previous example, let $G=\mathbb Z^2\ast \mathbb Z^2$ and let $H$ be the left hand side factor.  While $H$ is a qi-subcharacteristic subgroup by \Cref{cor:BDM-RelHyp}, the collection $\{H\}$ is not qi-characteristic. The second and third conditions of the \Cref{thm:qi-characteristic-groups2} hold, but the first does not. Specifically, a  quasi-isometry that flips the two factors sends $H$ to a space that is at infinite Hausdorff distance of any of its left cosets.

\subsection{Papasoglu's example}

 Consider an amalgamated product $G=\Z^3 \ast_\Z \Z^3$, where $\Z$ corresponds to a maximal infinite cyclic subgroup in both factors, and let $\mathcal{P}$ be the collection consisting of the two $\Z^3$ factors.  The work of Papasoglu~\cite[Theorem 7.1]{Pa05} implies that every $(L,C)$-quasi-isometry of $q\colon G\to G$ is $(L,C,M_q)$-quasi-isometry of pairs $q\colon (G, \mathcal{P}) \to (G, \mathcal{P})$ for some constant $M_q$. To show that $\mathcal{P}$ is qi-characteristic we need to show that $M_q$ can be chosen so that it depends only of $L$ and $C$, and not on $q$. \textbf{We do not know whether the constant $M_q$ can be chosen so that it only depends on $L$ and $C$}. Provided that is true, $\calP$ would be another example of a qi-characteristic collection.


\subsection{Baumslag-Solitar groups: a non-example}

Let $n\geq 2$ and consider the Baumslag-Solitar group $BS(1,n)=\langle a,t| tat^{-1}=a^n\rangle$.
The distorted cyclic subgroup $\langle a \rangle$ is not qi-characteristic since it has infinite index in its commensurator. 

The subgroup $\langle t\rangle$ does not form a qi-characteristic collection. We sketch the argument  using a construction that appears in the work of Farb and Mosher~\cite{FaMo98} on quasi-isometric rigidity of solvable Baumslag-Solitar groups. They use a particular metric on the Cayley complex $X_n$ of $BS(1,n)$ together with the projection $\pi\colon X_n\to T_n$ to the Bass-Serre tree. Let us recall a few properties: the inverse image $\pi^{-1}(L)$ 
of any \emph{coherently oriented  proper line $L$ of $T_n$} is an isometrically embedded hyperbolic plane $H$; all hyperbolic planes of $X_n$ arise in this way and can be simultaneously identified with the upper half plane model of $\mathbb{H}^2$ so that inverse image $\pi^{-1}(x)$ for $x\in L$ correspond to an  horocycle based at $\infty\in \partial \mathbb{H}^2$. In this way, the parabolic isometry $q\colon \mathbb{H}^2 \to \mathbb{H}^2$ given by  $z\mapsto z+1$ preserves horocycles based at $\infty$, and hence it induces an isometry $q\colon X_n \to X_n$ such that $\pi\circ q =\pi$. The isometry $q$ preserves each hyperbolic plane of $X_n$, and each of these planes corresponds to a unique left coset of $\langle t \rangle$ which can be identified with a particular vertical geodesic. Since any two hyperbolic planes of $X_n$ are at infinite Hausdorff distance, and any two distinct geodesics of $\mathbb{H}^2$ are at infinite Hausdorff distance, it follows that $q(\langle t\rangle)$ is at infinite Hausdorff distance of every left coset of $\langle t \rangle$.

\subsection{Hyperbolically embedded collections}
There is a characterisation of a subgroup $P\leq  G$ being hyperbolically embedded in terms of coned-off Cayley graphs \cite{MPR2021}.  Specifically, $P\hookrightarrow_h (G,T)$ if and only if $\hat\Gamma(G,P,T)$ is a hyperbolic graph which is fine at its cone vertices.  The authors have given sufficient conditions for a quasi-isometry of pairs to induce a quasi-isometry of coned-off Cayley graphs \cite[]{HuMaSa21} which leads to the following natural question.

\begin{question}
Suppose $\calp\hookrightarrow_h G$.  What conditions on a quasi-isometry of pairs $(G,\calp)\to(H,\calq)$ ensure $\calq\hookrightarrow_h H$?
\end{question}

There are already positive and negative results towards this question.  The first and second author have given technical conditions for a quasi-isometry of pairs to carry hyperbolically embedded collections in \cite{HMP202x}.  However, these conditions are often hard to apply since, in general, the generating set $T$ for which $P\hookrightarrow_h(G,T)$ is infinite and so constructing an infinite generating set for $H$ can be extremely tricky.  For an extended discussion the reader should consult loc. cit..

\begin{question}
Suppose $\calp\hookrightarrow_h G$.  What conditions ensure $\calp$ is qi-characteristic?
\end{question}

There is a result towards the previous question in the case of groups hyperbolic relative to NRH subgroups \cite{BDM09}.  However, the following example demonstrates the difficulty of the preceding question even in the case of finite extensions.

\begin{example}\cite{MiOs15,MiOs15c}\label{ex.MinasyanOsin}
Let $H=\langle a,b \rangle$ be the free group of rank two, let $G=\langle a,b,t \colon tat^{-1}=b,\quad t^2=e \rangle$, let $T=\{b,a,a^{-1},a^2,a^{-2},\ldots \}$ and $S=T\cup \{t\}$. The inclusion $\Gamma(H,T) \to \Gamma(G,S)$ is not a quasi-isometry. Indeed, in $G$ we have $ta^nt^{-1}=b^n$ and hence $\dist_{(G,S)}(e,b^n)=3$ but $\dist_{(H,T)}(e,b^n)=n$ for every $n$.  In particular, the map $\Gamma(H,T)\to \Gamma(H,T)$ given by $h\mapsto tht^{-1}$ is not a quasi-isometry, and hence the $G$-action on $H$ by conjugation is not an action by quasi-isometries.  
\end{example}

\section{Commensurated subgroups and coarse Poincar\'e duality}

The following section highlights recent work of Alexander Margolis which builds on work of many authors.  This first result shows that normal subgroups can be quasi-isometrically rigid even though they are not qi-characteristic.  To state the results in this section we say $q\colon(G,\calp)\to(H,\calq)$ is \emph{almost} a qi of pairs if $q$ is a quasi-isometry of pairs except that we do now require a uniform constant $M$ in the relation \eqref{eq:def-qi-relation}.

\begin{theorem} \emph{\cite[Theorem~5.10]{Margolis2021}}
Let $G$ be a finitely presented $\ZZ$-by-($\infty$ ended) group.  If $G'$ is any finitely generated group quasi-isometric to $G$, then it is also $\ZZ$-by-($\infty$ ended).  Moreover, any quasi-isometry $q\colon G\to G'$ induces a quasi-isometry $G/\ZZ\to G'/\ZZ$.
\end{theorem}

Note that the generalisation to $\ZZ^n$-by-($\infty$ ended) group fails if $n\geq2$.  Indeed, Leary--Minasyan have constructed groups quasi-isometric to $\ZZ^n\times F_2$ which do not have a normal $\ZZ^n$-subgroup \cite{LM2021}.  More examples quasi-isometric to certain RAAGs or right-angled buildings were constructed by the first author in \cite{Hughes2021} (see also \cite{Hughes2022}).

Coarse Poincar\'e duality groups were introduced by Kapovich and Kleiner as coarse analogues of Poincar\'e duality groups \cite{KapovichKleiner2005}.  The notion is invariant under quasi-isometries and has been studied by many authors \cite{MSW2003,Pa07,MSW11}.  A subgroup $H\leq G$ is \emph{commensurated} (or \emph{almost normal}) if every conjugate of $H$ in $G$ is commensurable with $H$, that is, for every $g\in G$ the indexes $|H\colon H\cap H^g|$ and $|H^g\colon H\cap H^g|$ are finite.

\begin{theorem}\emph{\cite[Theorem~1.5]{Margolis2021}}
If $G$ is a group of type $\mathsf{F}_{n+1}$ and $H$ is a commensurated coarse $\mathsf{PD}_n$ subgroup with $e(G/H)\geq3$, then every self quasi-isometry $q\colon G\to G$ extends to an almost quasi-isometry of pairs $q\colon (G,H)\to (G,H)$.
\end{theorem}

We defer the reader to \cite{Margolis2018} for the definitions of ``$(n-1)$-acyclic at infinity'' and ``coarsely $3$-separates''.  However, we include the following theorem (restated in the language of quasi-isometries of pairs) to illustrate that rigidity results for pairs about coarse homological properties are possible.

\begin{theorem}\emph{\cite[Theorem~1.3]{Margolis2018}}
Let $G$ be a group of type $\mathsf{FP}_{n+1}(\ZZ/2)$ that is $(n-1)$-acyclic at infinity over $\ZZ/2$.  Suppose $H\leq G$ is a coarse $\mathsf{PD}_n(\ZZ/2)$ group that coarsely $3$-separates $G$ and no infinite index subgroup of $H$ coarsely separates $G$.  Then every quasi-isometry $f\colon G\to G'$ extends to an almost quasi-isometry of pairs $f\colon(G,H)\to(G',H')$ such that $G'$ splits over $H'$.
\end{theorem}

Margolis also proves results about quasi-isometry of pairs preserving splittings regarding virtually polycylic groups \cite[Corollary~1.5]{Margolis2018} and fundamental groups of graphs of groups whose local groups satisfy appropriate finiteness and acyclicity at infinity conditions \cite[Corollary~1.4]{Margolis2018}.

\section{Sins of omission}
There are a number of topics we have not discussed in detail.  As we mentioned before, we have not discussed the theory of \emph{JSJ decompositions} and its connections to quasi-isometric rigidity.  We divert the readers attention to \cite{GuLe2017} for a comprehensive introduction and to \cite{DaniThomas2017,Barrett2018,TaTo2019,ReevesScottSwarup2020,ShWo2021,GardamKielakLogan2021} for recent developments and applications. We have not discussed \emph{coarse bundles} which have appeared in the work of many authors and serve as one of the main technical tools for proving rigidity results of extensions $N\cdot Q$.  The history of the construction can be traced through articles of Farb--Mosher \cite{FaMo99,FaMo2000,FaMo2002}, Whyte \cite{Whyte2001}, Mosher--Sageev--Whyte \cite{MSW2003}, Eskin--Fisher--Whyte \cite{EsFiWhy2012,EsFiWhy2013}, and most recently Margolis \cite{Margolis2021}.

\AtNextBibliography{\small}
\printbibliography

\end{document}